\documentclass[dvipsnames]{ifacconf}

\usepackage{algorithm}
\usepackage{algpseudocode}
\algdef{SE}[DOWHILE]{Do}{DoWhile}{\algorithmicdo}[1]{\algorithmicwhile\ #1}%
\usepackage{flushend}
\usepackage{amsmath}
\usepackage{IEEEtrantools}
\usepackage{mathtools}
\usepackage{./TLC}
\usepackage{blox}
\usepackage{physics}
\usepackage{graphicx}      
\usepackage{natbib}        

\usepackage{pgfplots}
\pgfplotsset{compat=1.16}
\usepgfplotslibrary{groupplots}
\usepackage{layouts}
\pgfplotscreateplotcyclelist{colors}{
        {CornflowerBlue},
        {BurntOrange},
        {RubineRed},
        {Emerald},
        {SpringGreen},
        {Goldenrod}
}

\begin{document}
\begin{frontmatter}

\title{Splitting algorithms and circuits analysis} 

\thanks[footnoteinfo]{The research leading to these results has received funding from the European
Research Council under the Advanced ERC Grant Agreement Switchlet n. 670645, and from
the Cambridge Philosophical Society.}

\author[Cambridge]{Thomas Chaffey} 
\author[KTH]{Amritam Das}
\author[Cambridge]{Rodolphe Sepulchre}

\address[Cambridge]{University of Cambridge, Department of Engineering, Trumpington Street,
Cambridge CB2 1PZ, {\tt\small \{tlc37, rs771\}@cam.ac.uk}.}
\address[KTH]{Division of Decision and Control Systems, School of Electrical Engineering and
Computer Science, KTH Royal Institute of Technology, 100 44 Stockholm, Sweden,
{\tt\small amritam@kth.se}.}

\begin{abstract}                
        The splitting algorithms of monotone operator theory find zeros of sums of
        relations.  This corresponds to solving series or 
        parallel one-port electrical circuits, or the negative feedback interconnection of two
        subsystems.  One-port circuits with series \emph{and} parallel
        interconnections, or block diagrams with multiple forward and return paths,
        give rise to current-voltage relations consisting of nested sums \emph{and
        inverses}.  In this
        extended abstract, we present new splitting algorithms specially suited to
        these structures, for interconnections of monotone and anti-monotone
        relations.
\end{abstract}


\end{frontmatter}
\section{Introduction}

The mathematical property of \emph{monotonicity} originated in the study of networks
of nonlinear
resistors \citep{Duffin1946, Zarantonello1960, Dolph1961, Minty1960, Minty1961,
Minty1961a}.  Monotonicity generalizes the concept of passivity from linear circuit
theory; loosely speaking, an element is monotone if it is passive with respect to
\emph{any} possible reference trajectory.  Following the influential paper of
\citet{Rockafellar1976}, monotone operator theory has grown to become a pillar of
large scale optimization theory \citep{Bauschke2011, Ryu2021a, Parikh2013,
Bertsekas2011}.  

Central to this theory are the family of
\emph{splitting algorithms}. These algorithms find zeros of sums of monotone operators, and allow
computation to be performed separately for each operator.  Recent work by the authors
has revisited the study of electrical networks using modern splitting algorithms
\citep{Chaffey2021a}.  The main idea is that finding a zero of the sum of two
operators is equivalent to solving the port behavior of their parallel (or series)
interconnection.  In turn, this is equivalent to solving the behavior of the negative
feedback interconnection of two elements.  This observation motivates the development
of splitting algorithms which match more general circuit architectures.  In
Section~\ref{sec:n}, we describe an algorithm which solves the behavior of arbitrary
series/parallel one-port circuits.

While splitting methods require each circuit element to be monotone, similar ideas
can be applied to \emph{mixed monotone} circuits, consisting of port interconnections
of monotone and anti-monotone elements.  This significantly expands the possible
types of circuit behavior, allowing, for example, relaxation oscillations
\citep{vanderPol1926} and
neuronal excitability \citep{FitzHugh1961}.  In \citep{Das2021}, the authors have adapted Difference of
Convex Programming \citep{Lipp2016, Yuille2003} to solve such behaviors.  In
Section~\ref{sec:difference}, we describe a new splitting algorithm which matches the
mixed feedback structure of oscillators such as the van der Pol and FitzHugh-Nagumo
models.

While classical splitting methods deal only with sums, the algorithms we describe
here deal with both sums and inverses - the two operations which constitute physical
port interconnections.  The algorithms described in this abstract form the basis for
a more general class of splitting algorithms, which correspond to arbitrary
interconnections of physical systems.

\section{Monotone and anti-monotone relations}
A Hilbert space $\mathcal{H}$ is a complete vector space equipped with an inner product,
$\bra{\cdot}\ket{\cdot}: \mathcal{H} \times \mathcal{H} \to \mathbb{C}$, and an 
induced norm $\norm{x} \coloneqq \sqrt{\bra{x}\ket{x}}$.  In this abstract, we will
treat general Hilbert spaces, although a common choice in practice is the space of
square-summable, discrete time signals on $[0, T]$, denoted $l_{2, T}$.

An \emph{operator} on $\mathcal{H}$, 
is a possibly
multi-valued map $R: \mathcal{H} \to \mathcal{H}$.  
The identity operator, which maps $u \in \mathcal{X}$ to itself, is denoted by $I$.
The domain of an operator $R$ is denoted $\dom{R}$.
The \emph{graph}, or \emph{relation}, of an operator, is the set $\{u, y\; | \; u \in
\dom{R}, y \in R(u)\} \subseteq \mathcal{H}\times\mathcal{H}$.  We use the notions of an operator and its relation
interchangeably, and denote them in the same way.  

The standard operations on functions can be extended to relations.  Let $R$ and $S$
be relations on an arbitrary Hilbert space $\mathcal{H}$.  Then:
\begin{IEEEeqnarray*}{rCl}
        S^{-1} &=& \{ (y, u) \; | \; y \in S(u) \}\\
        S + R &=& \{ (x, y + z) \; | \; (x, y) \in S, (x, z) \in R \}\\
        SR &=& \{ (x, z) \; | \; \exists\; y \text{ s.t. } (x, y) \in R, (y, z) \in S \}.
\end{IEEEeqnarray*}
Note that $S^{-1}$ always exists, but is not an inverse in the usual sense.  In
particular, in general $S^{-1}S \neq I$.

\begin{defn}
        A relation $S \subseteq \mathcal{H}\times\mathcal{H}$ is called \emph{monotone} if
\begin{IEEEeqnarray*}{rCl}
        \langle u_1 - u_2 | y_1 - y_2 \rangle \geq 0
\end{IEEEeqnarray*}
for any $(u_1, y_1), (u_2, y_2) \in S$.  
A monotone relation is called \emph{maximal} if it is not properly contained in any
other monotone relation.
\end{defn}

\begin{defn}
        A relation $S:\mathcal{H} \to \mathcal{H}$ is \emph{anti-monotone} if $-S$ is monotone.
\end{defn}

\section{Splitting two-element circuits}\label{sec:two}

There is a large body of literature on splitting algorithms, which solve problems of the form $0 \in M_1(u) +
M_2(u)$, where $M_1$ and $M_2$ are maximal monotone relations.  There is a direct analogy
with electrical circuits: if $M_1$ and $M_2$ are resistances, their series 
interconnection is given by the relation $v = M_1(i) + M_2(i)$; if $M_1$ and $M_2$
are conductances, their
parallel interconnection is given by $i = M_1(v) + M_2(v)$.  Given a current, 
the corresponding voltage across a parallel interconnection can be found using a
splitting algorithm, by solving $0 \in M_1(v) + M_2(v) - i$.
Here, we briefly describe two splitting algorithms -- the
forward/backward splitting, and the Douglas-Rachford splitting.  For the
convergence properties of these algorithms, we refer the reader to
\citep{Giselsson2019, Bauschke2011, Ryu2021a}.
Given an operator $S$ and a
scaling factor $\alpha$, the $\alpha$-resolvent of $S$ is defined to be the operator
\begin{IEEEeqnarray*}{rCl}
        \res_{\alpha S} \coloneqq (I + \alpha S)^{-1}.
\end{IEEEeqnarray*}
If $S$ is maximal monotone, $\res_{S}$ is single-valued \citep{Minty1961}.

\subsection{Forward/backward splitting}
The simplest splitting algorithm is the forward/backward splitting
\citep{Passty1979, Gabay1983, Tseng1988}.  Suppose $M_1$ and
$\res_{\alpha M_2}$ are single-valued.  Then:
\begin{IEEEeqnarray*}{lrCl}
       & 0 &\in& M_1(x) + M_2(x)\\
        \iff & 0 & \in & x - \alpha M_1(x) - (x + \alpha M_2(x))\\
        \iff &(I + \alpha M_2)x &\ni& (I - \alpha M_1)x \\
        \iff & x &=& \res_{\alpha M_2} (I - \alpha M_1) x.
\end{IEEEeqnarray*}
The fixed point iteration $x^{j+1} = \res_{\alpha M_2}(x^j - \alpha M_1(x^j))$ is the forward/backward splitting algorithm. 

\subsection{Douglas-Rachford splitting}

The reflected resolvent, or Cayley operator, is the operator
\begin{IEEEeqnarray*}{rCl}
        R_{\alpha S} \coloneqq 2\res_{\alpha S} - I.
\end{IEEEeqnarray*}
Given two operators $M_1$ and $M_2$, and a scaling factor $\alpha$, 
the Douglas-Rachford algorithm  \citep{Douglas1956, Lions1979} is the iteration
\begin{IEEEeqnarray*}{rCl}
z^{k + 1} &=& T(z^{k}),\\
x^k &=& \res_{\alpha M_2}(z^k),
\end{IEEEeqnarray*}
where $T$ is given by
\begin{IEEEeqnarray}{rCl}
T = \frac{1}{2}(I + R_{\alpha M_1} R_{\alpha M_2}).\label{eq:DR_operator}
\end{IEEEeqnarray}

\section{Splitting $n$-element circuits}\label{sec:n}

\begin{figure}[hb]
        \centering
        \includegraphics{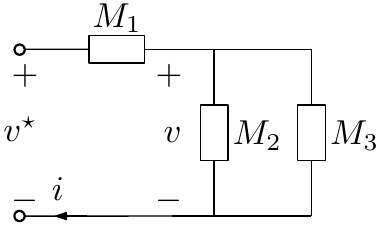}
        \caption{Three elements with one series
        interconnection and one parallel interconnection.}%
        \label{fig:three-circuit}
\end{figure}

If our circuit is composed of three elements, with one series interconnection and one
parallel interconnection (Figure~\ref{fig:three-circuit}), it has the form
$M = M_1 + (M_2 + M_3)^{-1}$. A naive approach to
solving the behavior of this circuit is to use a splitting algorithm such as the forward/backward algorithm,
with the resolvent step applied for $M_1$ and the forward step applied for $(M_2 +
M_3)^{-1}$.  Applying this forward step amounts to solving $v = (M_2 + M_3)^{-1}(i)$
for some $u$, which may be rewritten as $0 \in (M_2 + M_3)(v) - i$.  This can be
solved by again applying the forward/backward algorithm.

This naive procedure has poor complexity: for every forward/backward
step for $M_1 + (M_2 + M_3)^{-1}$, an \emph{entire} fixed point iteration has to
be computed for (an offset version of) $M_2 + M_3$.  In \citep{Chaffey2021a}, we propose an
alternative procedure for $n$-element circuits.  Here, we sketch this procedure on
the circuit of Figure~\ref{fig:three-circuit}. Rather than apply a forward step for the relation $(M_2 +
M_3)^{-1}$, we simply apply a \emph{single step} of the fixed point iteration needed
to compute this forward step, using the forward/backward algorithm.  Given $v^\star$, we want to solve 
$0 \in (M_1 + (M_2 + M_3)^{-1})(i) - v^\star$.  Assume that $M_3$, $\res_{\alpha_1 M_2}$ 
and $\res_{\alpha_2 M_1}$ are single-valued.  We then have:
\begin{IEEEeqnarray}{rCl}
        v^\star &\in& v + M_1(i)\label{eq:i_step}\\
        v &\in& (M_2 + M_3)^{-1}(i),\label{eq:v_step} 
\end{IEEEeqnarray}
where $v$ is the voltage over $M_2$, illustrated in
Figure~\ref{fig:three-circuit}.  Equation~\eqref{eq:i_step} gives
\begin{IEEEeqnarray*}{rCl}
        i + \alpha_2 M_1(i) &\ni& i - \alpha_2 v + \alpha_2 v^\star\\
        i &=& (I + \alpha_2 M_1)^{-1}(i - \alpha_2 v + \alpha_2 v^\star)\\
        i &=& \res_{\alpha_2 M_1}(i - \alpha_2 v + \alpha_2 v^\star).
\end{IEEEeqnarray*}
Equation~\eqref{eq:v_step} gives
\begin{IEEEeqnarray*}{rCl}
        i &\in& (M_2 + M_3)(v)\\
        v + \alpha_1 M_2(v) &\ni& v - \alpha_1 M_3 (v) + \alpha_1 i\\
        v &=& (I + \alpha_1 M_2)^{-1}(v - \alpha_1 M_3(v) + \alpha_1 i)\\
        v &=& \res_{\alpha_1 M_2}(v - \alpha_1 M_3(v) + \alpha_1 i).
\end{IEEEeqnarray*}
This shows that a fixed point of the iteration
\begin{IEEEeqnarray*}{rCl}
        v^{k + 1} &=& \res_{\alpha_1 M_2} (v^k - \alpha_1 M_3(v^k) + \alpha_1 i^k)\\
        i^{k + 1} &=& \res_{\alpha_2 M_1} (i^k - \alpha_2 v^{k+1} + \alpha_2 v^\star)
\end{IEEEeqnarray*}
is a solution to our original problem $0 \in (M_1 + (M_2 + M_3)^{-1})(i) - v^\star$.

\section{Splitting the difference}\label{sec:difference}

\begin{figure}[hb]
                \centering
                \includegraphics{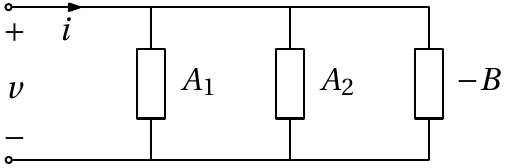}
                \caption{A parallel mixed monotone circuit, which is a prototype
                structure for systems such as the van der Pol and FitzHugh-Nagumo
        oscillators.}%
        \label{fig:mm_circuit}
        \end{figure}

A mixture of positive and negative feedback is a ubiquitous mechanism, in both biology
and engineering, for the generation of switches and oscillations \citep{Sepulchre2019,
Sepulchre2005, Stan2007, Stan2007a, Chua1987}.  Again adopting the analogy of
electrical circuits, such feedback systems can be thought of as the parallel
interconnection of three elements (Figure~\ref{fig:mm_circuit}) - the forward path and the negative feedback path,
which we assume to be monotone, and the positive feedback path, which we assume to be
anti-monotone.  Such a structure encompasses systems such as the van der Pol and FitzHugh-Nagumo
        oscillators.  For example, the van der Pol oscillator is given by $A_1(s) =
        (s^2 + 1)/s$, $A_2(v) = \mu v^3/3$ and $B(v) = \mu v$ (where $s$ is the
        Laplace variable) \citep{Das2021}.

Given the mixed monotone structure of Figure~\ref{fig:mm_circuit}, we can
find the steady state behavior of the system by solving a zero-finding
problem: $0 \in A_1(v) + A_2(v) - B(v) - i$.

The authors have explored methods to solve these problems using an adaptation of
Difference of Convex Programming in \citep{Das2021}.  The method involves iterating
the operator $(A_1 + A_2)^{-1} B$.  Computing $(A_1 + A_2)^{-1}$ at every iteration
is an expensive operation; in this section, we propose
the \emph{mixed monotone Douglas-Rachford algorithm} (Algorithm~\ref{alg:MMDR}), which replaces $(A_1 + A_2)^{-1}$ with a single step of the
Douglas-Rachford iteration needed to invert it. 

For operators $A_1$ and $A_2$ and step size $\alpha$, we define $T_\alpha (A_1, A_2)$ to be
the Douglas-Rachford operator:
\begin{IEEEeqnarray}{rCl}
T_\alpha (A_1, A_2) = \frac{1}{2}(I + R_{\alpha A_1} R_{\alpha A_2}).
\end{IEEEeqnarray}
Recall that $R_{\alpha S}$ denotes the Cayley operator $2\res_{\alpha S} - I$.

\begin{algorithm}[h!]
        \caption{Mixed-Monotone Douglas-Rachford}\label{alg:MMDR}
        \begin{algorithmic}[1]
                \State \textbf{Data:} Maximal monotone $A_1, A_2$. Monotone,
                single-valued $B$. Initial
                value $x_1$. Convergence tolerance $\epsilon > 0$.
                \State Define $A_1^j$ by $x \mapsto A_1(x) - y_j$ for all $j$.
                \State $j = 1$
                \Do 
                        \State Solve 
                                \begin{IEEEeqnarray*}{rCl}
                                        x_{j+1} &=& \res_{\alpha A_2}(z_{j})\\
                                        y_{j+1} &=& B(x_{j+1})\\
                                        z_{j+1} &=& T_\alpha (A_1^{j+1}, A_2a)(z_j).
                                \end{IEEEeqnarray*}
                        \State $j = j+1$.
                        \DoWhile{$|x_{j + 1} - x_j| > \epsilon$}
        \end{algorithmic}
        \label{alg:MMDR}
\end{algorithm}

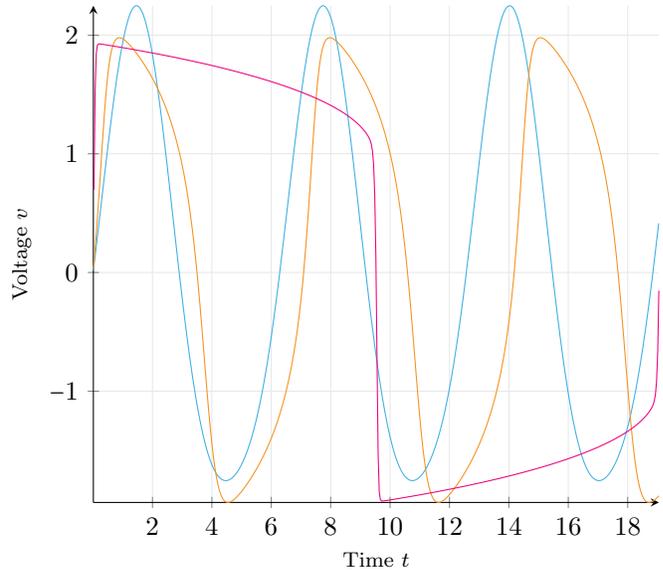
\begin{figure}[h]
        \centering
        \begin{tikzpicture}
                \begin{axis}
                        [
                        no markers,
                        name = ax1,
                        width=0.5\textwidth,
                        height=0.45\textwidth,
                        ticklabel style={/pgf/number format/fixed},
                        xlabel={\footnotesize Time $t$},
                        ylabel={\footnotesize Voltage $v$},
                        cycle list name=colors,
                        grid=both,
                        grid style={line width=.1pt, draw=Gray!20},
                        axis x line=bottom,
                        axis y line=left
                        ]
                        \addplot table [x = t, y = v, col sep = comma, mark =
                                none]{"./data/vdp_0.0002_reduced.csv"};
                                
                        \addplot table [x = t, y = v, col sep = comma, mark =
                                none]{"./data/vdp_1.5_reduced.csv"};

                        \addplot table [x = t, y = v, col sep = comma, mark =
                                none]{"./data/vdp_10_reduced.csv"};
               \end{axis}
        \end{tikzpicture}
        \caption{Steady-state solutions to the van der Pol oscillator for $\mu =
                0.0002$ (blue), $1.5$ (orange) and $10$ (red).
Algorithmic parameters are a step size of $\alpha = 0.05$, convergence tolerance of $\epsilon = 0.01$ and $5000$ time
steps.}%
\label{fig:vdp_solution}
\end{figure}

Note that a fixed point of this algorithm is a solution to $0 \in A_1(x) + A_2(x) -
B(x)$: we know, by convergence of the Douglas-Rachford algorithm, that $x$ is a
solution to $0 \in A_1^j(x) + A_2(x)$, which is equal to $A_1(x) + A_2(x) - B(x)$ at
a fixed point.  \citep[Thm. 4.1]{ChaffeyThesis} gives a convergence condition for
this algorithm.  Figure~\ref{fig:vdp_solution} shows steady-state solutions to the
van der Pol oscillator computed with Algorithm~\ref{alg:MMDR}.  The system is treated
as an interconnection of operators on $l_{2, T}$, the space of length $T$ periodic
signals.  For further details
of the implementation, the reader is referred to \citep[Example 4.3]{ChaffeyThesis}.

\bibliography{references}
\end{document}